\documentclass[12pt]{amsart}

\newtheorem{theorem}{Theorem}[section]

\theoremstyle{definition}

\theoremstyle{remark}

  {\end{list}}

\usepackage{amssymb} 
\usepackage{verbatim}

\def\Gal{{\mathop{\rm Gal}}}

\def\supp{{\mathop{\rm supp}}}

\def\dim{{\mathop{\rm dim}}}
\def\exp{{\mathop{\rm exp}}}
\def\vol{{\mathop{\rm vol}}}
\def\deg{{\mathop{\rm deg}}}

\def\Proj{{\mathop{\rm Proj}}}
\def\Spec{{\mathop{\rm Spec}}}

\def\GL{{\mathop{\rm GL}}}

\def\hdeg{{\widehat{\deg}}}
\def\Div{{\mathop{\rm div}}}

\def\<{\langle }
\def\>{\rangle }

\def\({(\!(}
\def\){)\!)}

\def\[{[\![}
\def\]{]\!]}

\def\BF1{{\mathbf 1}}

\def\fX{{\mathfrak{X}}}

\def\fp{{\mathfrak{p}}}

\def\AA{{\mathbb A}}
\def\CC{{\mathbb C}}
\def\QQ{{\mathbb Q}}
\def\EE{{\mathbb E}}

\def\RR{{\mathbb R}}
\def\ZZ{{\mathbb Z}}
\def\PP{{\mathbb P}}

\def\cO{{\mathcal O}}

\def\cK{{\mathcal K}}

\def\cL{{\mathcal L}}
\def\cF{{\mathcal F}}

\def\cY{{\mathcal Y}}
\def\cZ{{\mathcal Z}}

\def\cz{{Z}}

\def\tV{{\widetilde{V}}}

\def\fX{{\frak{X}}}

\def\zbar{{\overline{z}}}
\def\Hbar{{\overline{H}}}

\def\cLbar{{\overline{\cal L}}}

\def\Pbar{{\overline{P}}}

\def\sigmabar{{\overline{\sigma}}}
\def\cLbar{{\overline{\mathcal L}}}

\def\blacksquare{\Box} 

\DeclareMathSymbol{\varnothing} {\mathord}{AMSb}{"3F} 
 
\theoremstyle{definition} 
\theoremstyle{remark} 

\begin{document}

\title[Robin Formula]
{A Robin formula for the Fekete-Leja Transfinite Diameter}

\author{Robert Rumely}
\address{Robert Rumely\\ 
Department of Mathematics\\
University of Georgia\\
Athens, Georgia 30602\\
USA}
\email{rr@@math.uga.edu}

\date{September 21, 2005}
\thanks{Work supported in part by NSF grant DMS-0300784.}

\maketitle

The purpose of this note is to give a formula for the Fekete-Leja transfinite 
diameter on $\CC^N$, generalizing the classical Robin formula
\begin{equation*}
d_{\infty}(E) \ = \ e^{-V(E)}
\end{equation*}
for the usual transfinite diameter.

We will disengage this from a formula 
for the sectional capacity proved in arithmetic intersection theory
(\cite{CLR}, Theorem 1.1, p.233).  

\vskip .1 in
First recall the definition of the Fekete-Leja transfinite diameter for a 
compact set $E \subset \CC^N$ (see \cite{BC}, \cite{Z1}).  
Consider the set of monomials $z^{k} = z_1^{k_1} \cdots z_N^{k_N}$ 
in the polynomial ring $\CC[z] = \CC[z_1, \ldots, z_N]$.
Let $\Gamma(n)  \subset \CC[z]$ be the space of polynomials
of total degree at most $n$, 
and let $\cK(n) = \{k \in \ZZ^N : k_i \ge 0, k_1 + \cdots+ k_N \le n\}$
be the index set for the monomial basis of $\Gamma(n)$.
Put  $q_n = \#(\cK(n)) = \binom{n+N}{n}$.

Fixing $n$,  take $q_n$ independent vector variables 
$z_i = (z_{i1}, \ldots, z_{iN}) \in \CC^N$, $i = 1, \ldots, q_n$.
Let $k_1, \ldots, k_{q_n}$ be the indices in $\cK(n)$.
The Vandermonde determinant
\begin{equation*}
Q_n(z_1, \ldots, z_{q_n}) \ := \ \det(z_i^{k_j})_{i,j = 1}^{q_n}
\end{equation*}
is a homogeneous polynomial in the $z_{ij}$ of total degree
$T_n \ = N \cdot \binom{n+N}{N+1}$.  For each $n$, put
\begin{equation*}
d_n(E) \ = \
\max_{z_1, \ldots, z_{q_n} \in E} |Q_n(z_1, \ldots, z_{q_n})|^{1/T_n} \ .
\end{equation*}
The Fekete-Leja transfinite diameter is defined by 
\begin{equation*}
d_{\infty}(E) \ = \ \lim_{n \rightarrow \infty} d_n(E) \ .
\end{equation*}
The existence of the limit is due to Zaharjuta (\cite{Z1}).

Henceforth, we will assume that $d_{\infty}(E) > 0$.  
This is equivalent to $E$ being non-pluripolar (\cite{LT}).

\vskip .1 in

For $f \in \CC[z_1, \ldots, z_N]$ write $\|f\|_E = \sup_{z \in E} |f(z)|$.
The Green's function $G^*(z,E)$ is the
upper semicontinuous regularization of the Siciak extremal function
\begin{equation*}
G(z,E) \ := \ 
\lim_{n \rightarrow \infty} \max_{\substack{f \in \Gamma(n)  \\ \|f\|_E}}
            \frac{1}{n} \log(|f(z)|) \ .
\end{equation*}
Since $E$ is not pluripolar,
$G^*(z,E)$ is finite for all $z \in \CC^N$ and is plurisubharmonic
(\cite{Kl}, \cite{LT}, \cite{Z2}).
Write $dd^c = \frac{i}{4 \pi} \partial \overline{\partial}$ 
and let
\begin{equation} \label{FJ1}
\omega \ = \ dd^c (2G^*(z,E))
\end{equation}
be the associated positive $(1,1)$-current. (The factor $2$ is needed  
for compatibility with the Poincar\'e-Lelong formula).  

Let $\cz_0, \ldots, \cz_N$ be homogeneous coordinates on $\PP^N(\CC)$.
For $k=0, \ldots, N$ write $H_k$ for the hyperplane $\{ \cz_k = 0\}$  
and let $\AA^N_k = \PP^N(\CC) \backslash H_k$ 
be the corresponding affine patch.
Identify $\CC^N$ with $\AA^N_0 \subset \PP^N(\CC)$
via the embedding $(z_1, \ldots, z_N) \hookrightarrow (1 : z_1 : \cdots : z_N)$.
If $\|z\| = (|z_1|^2+ \cdots + |z_N|^2)^{1/2}$, then 
$|G^*(z,E) - \max(0, \log(\|z\|)|$ is uniformly bounded on $\CC^N$
(for boundedness from above, see \cite{Kl}, Corollary 5.2.2, p.193;  
for boundedness from below, note that since $E$ is compact, 
it is contained in a ball 
$B(0,R) = \{z \in \CC^N : \|z\| \le R\}$ for some large $R$, 
so $G^*(z,E) \ge G^*(z,B(0,R)) = \max(0,\log(\|z\|/R)$)).
It follows that for each $1 \le k \le N$ the function $g_k(z,E)$ defined on
$\AA^N_k \backslash H_0$ by
\begin{equation*}
g_k(z,E) \ = \ G^*(z,E) - \log(|z_k|)
\end{equation*}
extends uniquely to a plurisubharmonic function on $\AA^N_k$. 
If $k=0$, write $g_0(z,E) = G^*(z,E)$.  Since $\log(|\cz_i/\cz_j|)$
is pluriharmonic on $\AA^N_i \cap \AA^N_j$,
the currents $\omega_k := dd^c(2g_k(z,E))$ cohere to give 
a positive $(1,1)$-current on $\PP^N(\CC)$ extending $\omega$.  
We will denote it by $\omega$ as well.

Let $\omega^k = \omega \wedge \cdots \wedge \omega$ be the $k$-fold
exterior product of $\omega$ with itself, and put $\omega^0 = 1$.
Write $Y_k$ for the space $H_0 \cap \cdots \cap H_{k-1}$ and let  
$U_k = Y_k \backslash H_k$, so $U_k(\CC) \cong \CC^{N-k}$.  
Define the {\it iterated Robin constant} by
\begin{eqnarray}
\tV(E) & = & \frac{1}{N} \sum_{k=1}^N
    \int_{U_k(\CC)} g_k(z,E) \, \omega^{N-k} \label{F2} \\
    & = & \frac{1}{N} \big( \int_{Z_0=0,Z_1 \ne 0} 
               g_1(z,E) \, \omega^{N-1} \notag \\
    &  & \qquad \qquad           
      + \int_{Z_0=Z_1=0, Z_2 \ne 0} g_2(z,E) \, \omega^{N-2} \\
    & & \qquad \qquad \qquad \qquad  
                    + \cdots + g_N((0:\cdots:0:1),E) \big) \notag \ .
\end{eqnarray}
Note that when $N=1$, the sum consists of a single term and
reduces to the usual Robin constant, since $\infty = (0:1) \in \PP^1(\CC)$, 
and $g_1(\infty,E) = \lim_{z \rightarrow \infty} G^*(z,E) - \log(|z|) = V(E)$.  

\begin{theorem} \label{Thm1}
If $E \subset \CC^N$ is compact and not pluripolar, then
\begin{equation} \label{F3A}
d_{\infty}(E) = e^{-\tV(E)} \ .
\end{equation}
\end{theorem}

\vskip .1 in
Before giving the proof we will need some facts 
from arithmetic intersection theory.  Let $K$ be a number field, and let
$\cO_K$ be the ring of integers of $K$.  
Let $X/K$ be a smooth, connected projective variety of dimension $N$,
and write $K(X)$ for the field of $K$-rational functions on $X$.  

\vskip .1 in
  
The sectional capacity $S_{\gamma}(\EE,D)$ 
is a measure of size for an adelic set $\EE$ on $X$
relative to an ample divisor $D$. 
It was first proposed by Chinburg (\cite{Ch}),  
and its existence was shown in (\cite{RLV}), using methods from (\cite{Z1}).  
The name `sectional capacity' 
refers to asymptotics of volumes of spaces of sections of 
$\cO_X(nD)$ as $n \rightarrow \infty$.   

Let $D$ be an effective, ample $K$-rational Cartier divisor on $X$. 

For each place $v$ of $K$, let $K_v$ be the completion of $K$ at $v$, and let
$\CC_v$ be the minimal complete algebraically closed field containing $K_v$.
Let $|x|_v$ be the absolute value on $\CC_v$ extending the canonical
absolute value on $K_v$ given by the modulus of additive Haar measure.  
Without loss, we can assume that $X$ is embedded in $\PP^M$ for a suitable $M$.
There is a natural distance function $d_v(x,y)$ on $\PP^M(\CC_v)$:  
the chordal metric associated to the Fubini-Study metric, 
if $v$ is archimedean;  the $v$-adic spherical metric, 
if $v$ is nonarchimendean (see \cite{R1}, \S1.1).  

For each $v$, let $E_v \subset X(\CC_v)$ be a nonempty set,
and write $\EE = \prod_v E_v$;  we will call $\EE$ an {\it adelic set}.
We will assume
that the $E_v$ and $D$ satisfy the following `Standard Hypotheses':

\vskip .05 in  
    (1) Each $E_v$ is bounded away from $\supp(D)(\CC_v)$ under $d_v(x,y)$ 
and is stable under the group of continuous automorphisms $\Gal^c(\CC_v/K_v)$.  

    (2) For all but finitely many $v$, $E_v$ and $\supp(D)(\CC_v)$  specialize
to disjoint sets$\pmod{v}$;  equivalently, for all but finitely many $v$, 
$d_v(x,y) = 1$ for all $x \in E_v$ and all $y \in \supp(D)(\CC_v)$.  
\vskip .05 in 

\noindent{Note} that if $v$ is archimedean and $K_v \cong \RR$, 
then $\Gal^c(\CC_v/K_v) = \{1,\tau\}$ where $\tau$ is complex conjugation;  
if $K_v \cong \CC$, then $\Gal^c(\CC_v/K_v)$ is trivial.  

For each integer $n \ge 0$, put $\Gamma(nD) = H^0(X,\cO_X(nD))$,  
a finite dimensional vector space over $K$.  
Dehomogenizing at $D$, identify $\Gamma(nD)$ with 
$\{ f \in K(X) : \Div(f) + nD \ge 0 \}$.  
For each $v$, consider the space of
$K_v$-rational functions on $X$ with polar divisor at most $nD$,
$\Gamma_v(nD) = K_v \otimes_K \Gamma(nD) 
= \{f \in K_v(X) : \Div(f) + nD \ge 0 \}$.  For $f \in K_v(X)$,
write $\|f\|_{E_v} = \sup_{x \in E_v} |f(x)|_v$ and let 
\begin{equation*}
\cF(E_v,nD) \ = \ \{f \in \Gamma_v(nD) : \|f\|_{E_v} \le 1\} \ .
\end{equation*}
Let $K_{\AA} = \{(x_v) \in \prod_v K_v : \text{$|x_v|_v \le 1$ for all but 
finitely many $v$}\}$ be the adele ring of $K$.  
Identifying $K_{\AA} \otimes_K \Gamma(nD)$ with a subset of 
$(\prod_v K_v) \otimes_K \Gamma(nD) \cong \prod_v (K_v \otimes_K \Gamma(nD))$,
introduce the `adelic unit ball' 
\begin{equation*}
\cF(\EE,nD) \ = \ 
      (K_{\AA} \otimes_K \Gamma(nD)) \, \bigcap \, (\prod_v \cF(E_v,nD)) \ .
\end{equation*}
Fix an additive Haar measure $\vol_{\AA}$ on $K_{\AA}$.  
By transport of structure using a $K$-basis for $\Gamma(nD)$,
$\vol_{\AA}$ induces a Haar measure on $K_{\AA} \otimes_K \Gamma(nD)$.   
The Product Formula shows this measure is independent
of the choice of basis;  it too will be denoted $\vol_{\AA}$.
The sectional capacity $S_{\gamma}(\EE,D)$ is defined by 
\begin{equation*}
-\log(S_{\gamma}(\EE,D)) \ = \ 
 \lim_{n \rightarrow \infty} \frac{(N+1)!}{n^{N+1}} 
         \vol_{\AA}(\cF(\EE,nD)) 
\end{equation*}
(see \cite{RLV}, pp.23-24).
The limit is independent of the choice of Haar measure on $K_{\AA}$.  
Its existence is part of (\cite{RLV}, Theorem C, p.8).  

\vskip .1 in
 
We will now discuss Arakelov theory.  
Let $\fX$ be a model of $X$:  an integral projective scheme of 
Krull dimension $N+1$, flat and proper over $\Spec(\cO_K)$, 
whose generic fibre is $K$-isomorphic to $X$.

Put $X_{\CC} = X \times_\QQ \Spec(\CC)$;  
then $X_{\CC}(\CC) = \bigsqcup_{\sigma} X_{\sigma}(\CC)$  
is a complex manifold with a component corresponding to each embedding 
$\sigma : K \hookrightarrow \CC$.  
Write $\tau$ for complex conjugation on $\CC/\RR$;
there is a natural operation of $\tau$ on $X_{\CC}(\CC)$ 
coming from its action on the base.  If $\sigma$ is a real embedding, then 
$\tau$ acts as complex conjugation on $X_{\sigma}(\CC)$;  
if $\{\sigma, \sigmabar\}$ is a pair of conjugate complex embeddings,
then $\tau$ interchanges the components $X_{\sigma}$ and $X_{\sigmabar}$,
mapping $P \in X_{\sigma}(\CC)$ to $\Pbar \in X_{\sigmabar}(\CC)$, 
where $\Pbar$ is the point whose coordinates are the complex conjugates of 
those of $P$.    

Let $L$ be a line bundle on $X$.  
A {\it metrized line bundle}
$\cLbar = (\cL,\[ \cdot \])$ on $\fX$, extending $L$, 
is a pair consisting of a locally free sheaf $\cL$ of rank $1$ on $\fX$ 
which induces $L$ on $X$,  
and a smooth positive Hermitian metric $\[ \cdot \]$ 
on the fibres of $L$ over $X_{\CC}(\CC)$,  
which is invariant under the action of $\tau$.
Let $c_1(\cLbar)$ be the first Chern class of $\cLbar$,
the smooth $(1,1)$-form on $X_{\CC}(\CC)$ 
defined in a neighborhood of any point $z_0 \in X_{\CC}(\CC)$ 
by $dd^c(-2\log(\[s(z)\]))$, where $s$ is a meromorphic section of $L$ which
is defined and does not vanish at $z_0$.  

For each cycle $\cY$ on $\fX$, Bost, Gillet and Soul\'e (\cite{BGS}) define
the ``height'' $h_{\cLbar}(\cY)$ of $\cY$ relative to $\cLbar$ 
to be the self-intersection number
$\hdeg(\hat{c}_1(\cLbar)^{\dim(\cY)}|_{\cY})$
in the arithmetic intersection theory Gillet and Soul\'e 
developed in (\cite{GS2}, \cite{GS3}, \cite{GS4}).
The height is additive in $\cY$.
There is a recursive formula for $h_{\cLbar}(\cY)$
(see \cite{BGS}, Proposition 3.2.1(iv), p.949,
and the remarks after it dealing with the non-regular case):  
when $\cY$ is irreducible,
let $s \ne 0$ be a section of $\cL|_{\cY}$, 
and let $\Div_{\cY}(s)$ be its divisor.  
Write $Y$ for the generic fibre $\cY_K$.  Then
\begin{equation} \label{RecForm}
h_{\cLbar}(\cY) \ = \ \int_{Y_{\CC}(\CC)} -\log(\[s(z)\]) \, c_1(\cLbar)^{\dim(Y)}
                     + h_{\cLbar}(\Div_{\cY}(s)) \ . 
\end{equation}    
Eventually the recursion abuts at a $0$-cycle $\cZ = \sum_{\fp} n_{\fp} \cdot \fp$, 
a finite sum of closed points, and then 
$h_{\cLbar}(\cZ) = \sum n_{\fp}  \log(N\fp)$ where $N\fp$ 
is the order of the residue field at $\fp$ (\cite{BGS}, formula 3.1.4, p.946).
 
It is customary to write $\cLbar^{N+1}$ for $h_{\cLbar}(\fX)$.  
If $\cLbar$ is replaced by $\cLbar^{\otimes m}$, then 
$(\cLbar^{\otimes m})^{N+1} = m^{N+1} \cdot \cLbar^{N+1}$  
(\cite{BGS}, Proposition 3.2.1(i)).  This leads
to the notion of a {\it fractional metrized line bundle}:
if $n > 0$ is an integer, and $\cLbar$ is a metrized line bundle 
inducing $L^{\otimes n}$ on $X$, then we call the formal object 
$\frac{1}{n}\cLbar$ a fractional metrized line bundle,
and define $(\frac{1}{n}\cLbar)^{N+1} = n^{-(N+1)} \cdot \cLbar^{N+1}$.  
  
\vskip .1 in

Chinburg, Lau and Rumely (\cite{CLR}) expressed the sectional capacity 
as a limit of self-intersection numbers of fractional metrized line bundles.
Given an adelic set $\EE$ and an effective ample divisor $D$ on $X$ 
satisfying the Standard Hypotheses, 
they constructed a sequence of models $\fX_n$,
and metrized line bundles $\cLbar_n$ 
on $\fX_n$ extending $L^{\otimes n} = \cO_X(nD)$ on $X$, such that
\begin{equation} \label{FKK0}
-\log(S_{\gamma}(\EE,D)) 
\ = \ \lim_{n \rightarrow \infty} (\frac{1}{n}\cLbar_n)^{N+1} \ :
\end{equation}
see (\cite{CLR}, Theorem 1.1, p.233). 
 
The models $\fX_n$ are defined using the nonarchimedean part of $\EE$.
For each $n$, put $S_n = 
\{ f \in \Gamma(nD) : \text{$\|f\|_{E_v} \le 1$ for all nonarchimedean $v$}\}$,
and let $\cO_K[S_n]$ be the graded 
$\cO_K$-algebra generated in degree $1$ by $S_n$.
Then $\fX_n = \Proj(\cO_K[S_n])$.  Let $\cL_n = \cO_{\fX_n}(1)$; 
then $\cL_n$ induces $L^{\otimes n}$ on $X$. 

The metrics $\[ \cdot \]_n$ are constructed using the 
archimedean part of $\EE$.
For each archimedean place $v$ of $K$, fix an isomorphism $\CC_v \cong \CC$,  
and choose a sequence
of sets $E_{v,1} \supseteq E_{v,2} \supseteq \cdots $  containing $E_v$
for which the extremal functions  
\begin{equation*}
G(z,E_{v,n},D) \ := \ \lim_{m \rightarrow \infty}
   \sup_{\substack{ f \in \Gamma(mD) \\ \|f\|_{E_{v,n}} \le 1 }}
           \frac{1}{m} \log(|f(z)|)
\end{equation*}
are continuous and increase monotonically to $G(z,E_v,D)$.
(The existence of such sets $E_{v,n}$ follows from the proof of
(\cite{CLR}, Lemma 1.2, p.234).)

Since $G(z,E_{v,n},D)$ is continuous, it is plurisubharmonic.
By a theorem of Richburg (\cite{Ri}, Satz 4.7) 
there is a smooth plurisubharmonic function $G_{v,n}(z)$ such that 
\begin{equation*}
G(z,E_{v,n},D) - \frac{1}{n} \ \le \ G_{v,n}(z) 
             \ \le \ G(z,E_{v,n},D)-\frac{1}{n+1}  
\end{equation*}
for all $z \in X(\CC_v) \backslash \supp(D)(\CC_v)$.  Hence
\begin{equation*}
G(z,E_v,D) \ = \ \lim_{n \rightarrow \infty} G_{v,n}(z) 
\end{equation*}
as an increasing limit.  

By (\cite{CLR}, Theorem 2.13, p.253), the smoothings
$G_{v,n}(z)$ can be chosen so that for each $z_0 \in \supp(D)(\CC_v)$, 
if $s$ is a local equation for $D$ at $z_0$,  
then $G_{v,n}(z)+\log(|s(z)|)$ extends to a smooth
plurisubharmonic function in a neighborhood of $z_0$.  
If $K_v \cong \RR$, then since $E_v$ has been assumed to be stable
under complex conjugation, the $G_{v,n}(z)$ can be chosen to be invariant
under complex conjugation as well.  

Each embedding $\sigma : K \hookrightarrow \CC$ determines a place $v$ of $K$.
If $\sigma$ is a real embedding, put $G_{\sigma,n}(z) = G_{v,n}(z)$.  
If $\sigma$ is a complex embedding, 
then precisely one of $\sigma$ and its complex conjugate
$\sigmabar$ induces the chosen isomorphism $\CC_v \cong \CC$;  if it is $\sigma$, 
put $G_{\sigma,n}(z) = G_{v,n}(z)$, 
if not, put $G_{\sigma,n}(z) = G_{v,n}(\zbar)$.  

The metric
$\[ \cdot \]_n$  is defined 
by requiring that for the tautological section $1$ of $L^{\otimes n}$,
if $z \in X_{\sigma}(\CC) \subset X_{\CC}(\CC)$ then   
\begin{equation*}
\[1(z)\]_n \ = \ \exp(-n G_{\sigma,n}(z)) \ .
\end{equation*}  
By construction, $\[ \cdot \]_n$ is invariant under $\tau$.

\vskip .1 in
\noindent{\bf Proof of Theorem \ref{Thm1}.}
Let $E \subset \CC^N$ be a compact, non-pluripolar set.   
We will embed it as the archimedean component of an adelic set,
and apply the machinery above.

Take $K = \QQ(\sqrt{-1})$ as the ground field, so $\cO_K$ is the ring
of Gaussian integers.  Write $\PP^N_{\cO_K}$ for $\PP^N/\Spec(\cO_K)$, 
and write $\PP^N_{K}$ for its generic fibre.
Let $X/K$ be $\PP^N_{K}$, and let $D$ be defined by $\{Z_0 = 0\}$.
 
Put $\AA^N = \PP^N_{K} \backslash H_0$, 
and define $\EE = \prod_v E_v \subset \prod_v \AA^N(\CC_v)$ as follows.
For each place $v$, identify $\AA^N(\CC_v)$ with $\CC_v^N$. 
There is one archimedean place $v_\infty$ for $K$; 
fix an isomorphism $K_{v_\infty} \cong \CC$ 
and put $E_{v_\infty} = E \subset \CC^N$.
Condition (1) of the Standard Hypotheses holds trivially.
For each nonarchimedean $v$, put 
$E_v = B(0,1) = \{(z_1, \cdots, z_N) \in \CC_v^N : \max(|z_i|_v) \le 1\}$.  
Since $\Gal^c(\CC_v/K_v)$ preserves $|x|_v$, 
again condition (1) of the Standard Hypotheses holds.
The sets $E_v$ and $H_0(\CC_v)$ specialize
to disjoint sets$\pmod{v}$ for all nonarchimedean $v$,
so condition (2) of the Standard Hypotheses holds as well.

By (\cite{RL}, Theorem 3.1, p.551) the sectional capacity $S_{\gamma}(\EE,D)$
can be decomposed as a product of `local sectional capacities'
\begin{equation*}
S_{\gamma}(\EE,D) \ = \ \prod_v S_{\gamma}(E_v,D) \ .
\end{equation*}
Here the $S_{\gamma}(E_v,D)$ depend on the choice of an ordered basis
for the graded ring $\oplus_{n=0}^{\infty} \Gamma(nD)$: 
we take this to be the monomial basis,
equipped with the lexicographic order graded by the degree.
In this situation, $S_{\gamma}(E_{v_\infty},D) \ = \ d_{\infty}(E)^{2N}$.
This follows from (\cite{RL}, Theorems 2.3 and 2.6), 
combined with the discussion on (\cite{RL}, p.557).  
(In (\cite{RL}, p.557), the formula  $S_{\gamma}(E_{\infty},D) = d_{\infty}(E)^N$
is given when $\RR$ is the ground field.  Here, since $K_{v_\infty} \cong \CC$,
the normalized absolute value $|x|_{v_\infty} = |x|^2$ 
used in computing $S_{\gamma}(E_{v_\infty},D)$ 
is the square of the usual absolute value;  
this is the source of the $2$ in the exponent.) 
Furthermore, $S_{\gamma}(E_v,D) = 1$ for each nonarchimedean $v$
by (\cite{RL}, Example 4.1, p.555).   Hence
\begin{equation} \label{FGG1} 
S_{\gamma}(\EE,D) \ = \ d_{\infty}(E)^{2N} \ .
\end{equation}

Let $\{\frac{1}{n}\cLbar_n\}$ be the sequence of fractional 
metrized line bundles
constructed in (\cite{CLR}, Theorem 1.1) as described above.  
Since $E_v$ is the `trivial set' $B(0,1)$ for each nonarchimedean $v$,
the Maximum Modulus Principle of nonarchimedean analysis shows that 
\begin{equation*}
S_n \ = \ 
\oplus_{k_0 + \cdots + k_N = n} \  
      \cO_K \cdot Z_0^{k_0} \cdots Z_N^{k_N} \ .
\end{equation*}
It follows that $\fX_n \cong \PP^N_{\cO_K}$, 
and $\cL_n \cong \cO_{\PP^N_{\cO_K}}(n)$.  
Let $\cLbar_n^{\prime}$ be $\cO_{\PP^N_{\cO_K}}(1)$, equipped with the metric
$\[ \cdot \]_n^{\prime}$ defined by $\[1(z)\]_n^{\prime} = \exp(-G_n(z))$,
where $G_n(z)$ is metric constructed as above by smoothing $G(z,E,D)$.  
Then $\cLbar_n \cong (\cLbar_n^{\prime})^{\otimes n}$, 
so $(\frac{1}{n}\cLbar_n)^{N+1} = (\cLbar_n^{\prime})^{N+1}$
and without loss we can replace $\frac{1}{n}\cLbar_n$ by $\cLbar_n^{\prime}$.

We will now compute the intersection product $(\cLbar_n^{\prime})^{N+1}$
using Bost-Gillet-Soul\'e's recursive formula.  
There are two embeddings of $K$ into $\CC$,
both of which correspond to the place $v_{\infty}$, 
so $(\PP^N_K)_{\CC}(\CC)$ has two 
components which are interchanged by $\tau$.  
Write $\PP^N(\CC)$ for the one corresponding
to our chosen isomorphism $K_{v_\infty} \cong \CC$,
and $(\PP^N)_{\tau}(\CC)$ for the other.  Since they are isomorphic, 
and the metrics $\[ \cdot \]_n$ are $\tau$-invariant, the integrals over the two
components in the archimedean part of the intersection product are the same.  
Therefore, in what follows, we will compute the integrals over $\PP^N(\CC)$ 
and double the answer.  For any $K$-rational subvariety $Y \subset \PP^N_K$, 
write $Y(\CC)$ for the part of $Y_\CC(\CC)$ in the chosen component $\PP^N(\CC)$. 
  
Write $\AA^N_k$ for $\PP^N_K \backslash H_k$, where $H_k = \{Z_k = 0\}$.
Let $\Hbar_k$ be the Zariski closure of $H_k$ in $\PP^N_{\cO_K}$.
Then $\Hbar_0, \cdots, \Hbar_N$ meet transversely on $\PP^N_{\cO_K}$.
Write $g_{n,0}(x) = G_n(x)$ on $(\AA^N_0)(\CC)$ 
and for each $k = 1, \ldots, N$ let $g_{n,k}(x)$ be the
natural extension of $G_n(x)- \log(|z_k(x)|)$ to a plurisubharmonic function
on $\AA^N_k(\CC)$;  here $z_k = Z_k/Z_0$.  
As before, the $(1,1)$-forms $dd^c(2g_{n,k})$
glue to give a well-defined $(1,1)$-form $\omega_n$ on $\PP^N(\CC)$.  

Since $L = \cO_{\PP^N_K}(D)$ where $D = \{Z_0 = 0\}$, 
the canonical section `$1$' of $L$ is $Z_0$.    
This means that for each $x \in \AA^N_0(\CC)$ 
\begin{equation*}
-\log(\[Z_0(x)\]_n^{\prime}) \ = \ G_n(x) \ = \ g_{n,0}(x) \ . 
\end{equation*}
Similarly, on $\AA^N_k(\CC)$, 
\begin{eqnarray*}
-\log(\[Z_k(x)\]_n^{\prime}) 
& = & -\log(\[(Z_k/Z_0)(x) \cdot Z_0(x)\]_n^{\prime}) \\
                    & = & -\log(|z_k(x)|) + G_n(x) \ = \ g_{n,k}(x) \ .  
\end{eqnarray*}

Put $\cY_0 = \fX_n = \PP^N_{\ZZ}$ and take $s = Z_0$.  
Let $Y_0 = \PP^N_{\QQ}$ be the generic fibre of $\cY_0$ and 
put $\cY_1 = \Div(Z_0) = \Hbar_0$. 
By Bost-Gillet-Soul\'e's formula (\ref{RecForm}) and the remarks above, 
\begin{equation*}
(\cLbar_n^{\prime})^{N+1} \ = \ h_{\cLbar_n^{\prime}}(\fX_n) 
\ = \ 2 \int_{Y_0(\CC)} g_{n,0}(z) \, \omega^N 
\ + \  h_{\cLbar_n^{\prime}}(\cY_1) \ .
\end{equation*}                    
Inductively apply this formula 
to the sections $Z_1, \cdots, Z_N$, putting 
$\cY_{k+1} = \Div_{\cY_k}(Z_k|_{\ \cY_k}) 
= \Hbar_0 \cdot \, \ldots \, \cdot \Hbar_k$.
Note that $\Hbar_0 \cdot \, \ldots \, \cdot \Hbar_N = 0$, so the 
final abutment term vanishes.  It follows that 
\begin{equation} \label{FKK1}
(\cLbar_n^{\prime})^{N+1}  \ = \
2 \sum_{k=0}^N \int_{Y_k(\CC)} g_{n,k}(z) \, \omega_n^{N-k} \ .
\end{equation}

Since $\omega_n$ is smooth, the $(N-k-1)$-dimensional subspace $Y_{k+1}(\CC)$ 
of $Y_k(\CC)$ has measure $0$ under $\omega_n^{N-k}$.  
Put $U_k = Y_k \backslash Y_{k+1}$.  Then for each $k$ 
\begin{equation} \label{FKK2}
\int_{Y_k(\CC)} g_{n,k}(z)\, \omega_n^{N-k} 
\ = \ \int_{U_k(\CC)} g_{n,k}(z)\, \omega_n^{N-k} \ .
\end{equation}

Now let $n \rightarrow \infty$.  
For each $k$ the $g_{n,k}$ increase monotonically on $U_k(\CC)$ 
to $G(z,E)-\log(|z_k|)$, whose upper semicontinuous regularization
is $g_k(z,E)$.  Noting that $U_k(\CC) \cong \CC^{N-k}$,
it follows from (\cite{BT}, Theorem 7.4) that 
\begin{equation} \label{FKK3} 
\lim_{n \rightarrow \infty} \int_{U_k(\CC)} g_{n,k}(z)\, \omega_n^{N-k}
\ = \ \int_{U_k(\CC)} g_k(z,E)\, \omega^{N-k} \ .
\end{equation}
Combining (\ref{FKK0}), (\ref{FKK1}), (\ref{FKK2}), and (\ref{FKK3}) gives  
\begin{equation} \label{FGG2}
-\log(S_{\gamma}(\EE,D)) \ = \ 
2\sum_{k=0}^N \int_{U_k(\CC)} g_k(z,E) \, \omega^{N-k} \ .
\end{equation}

When $k = 0$,  $g_0(z,E)$ is the extremal 
plurisubharmonic function $G^*(z,E)$ on $U_0(\CC) \cong \CC^N$, 
and $\omega^N$ is the Bedford-Taylor measure, 
which is supported on $E$.
The set $\{z \in E : G^*(z,E) > 0\}$ is a negligible set, 
so it has measure $0$ under $\omega^N$.  Hence 
\begin{equation} \label{FGG3}
\int_{U_0(\CC)} g_0(z,E) \, \omega^N \ = \ 0 \ .
\end{equation}
Using (\ref{FGG1}), (\ref{FGG2}), and (\ref{FGG3}) we obtain 
\begin{equation*}
-2N \cdot \log(d_{\infty}(E)) \ = \ 2\sum_{k=1}^N \int_{U_k(\CC)} g_k(z,E) \, \omega^{N-k}
\end{equation*} 
which is equivalent to (\ref{F3A}).  \hfill $\blacksquare$

\vskip .1 in
\noindent{\bf Generalizations.}
It is tempting to assert that the formula in Theorem \ref{Thm1} 
gives a new  definition of the capacity.  However, we do not 
do so because the stated formula is only one of a class of 
formulas with the same property.  

In particular, the order in which the hyperplanes 
$H_k$ are intersected in the definition of $\tV(E)$ is immaterial.  
More generally, if $A = (a_{kj}) \in \GL_N(\CC)$, 
and if $Z_k^{\prime} = \sum a_{kj}Z_j$ for $k = 1,\ldots,N$
then the constant $\tV(E)$ can equally well be defined by 
\begin{equation} \label{FGL1} 
\tV(E) \ = \ \frac{1}{N} \sum_{k=1}^N 
       \int_{Z_0 = Z_1^{\prime}=\cdots=Z_{k-1}^{\prime}=0, Z_k^{\prime} \ne 0} 
                 g_k^{\prime}(z,E) \, \omega^{N-k} \ 
        + \frac{1}{N} \log(|\det(A)|) 
\end{equation}   
where $g_k^{\prime}(z,E) = G^*(z,E)-\log(|\sum a_{kj}z_j|)$ on $\CC^N$.    

Indeed, 
if we put $z^{\prime} = {}^t(z_1^{\prime}, \ldots, z_N^{\prime}) =  A(z)$, 
where $z_k^{\prime} = \sum a_{kj}z_j$ for $k =1, \ldots, N$,
this follows from She\u{i}nov's formula 
(see \cite{Shei}, or \cite{BC}, p.287)
\begin{equation*} 
d_{\infty}(A(E)) \ = \ |\det(A)|^{1/N} \cdot d_{\infty}(E) \ , 
\end{equation*} 
by applying Theorem \ref{Thm1} to the set $A(E)$, 
and noting that $G(z^{\prime},A(E)) = G(z,E)$, 
as follows easily from the definitions.

\end{document}